\documentclass[12pt]{article}
\usepackage{amsmath,amssymb,amsfonts,graphicx,wrapfig}

\newtheorem{theorem}{Theorem}[section]
\newtheorem{lemma}[theorem]{Lemma}
\newtheorem{claim}[theorem]{Claim}

\newcommand{\cl}{{\mathrm{cl}}}
\newcommand{\ind}{{\mathrm{ind}}}

\newcommand{\qd}{{\mathrm{q}}}

\newcommand{\const}{{\mathrm{const}}}

\begin{document}

\begin{center}
{\Large
Logical laws for short existential monadic second order sentences about graphs
\\}

\vspace{0.2cm}

{\large
M.E. Zhukovskii~\footnote{Moscow Institute of Physics and Technology (National Research University), laboratory of advanced combinatorics and network applications, 9 Institutskiy per., Dolgoprodny, Moscow Region, 141701, Russian Federation; Adyghe State University, Caucasus mathematical center, ul. Pervomayskaya, 208, Maykop, Republic of Adygea, 385000, Russian Federation; The Russian Presidential Academy of National Economy and Public Administration, Prospect Vernadskogo, 84, bldg 2, Moscow, 119571, Russian Federation.}
}
\end{center}

\vspace{0.5cm}

{\small
\begin{center}
Abstract
\end{center}
\vspace{0.1cm}

In 2001, J.-M. Le Bars proved that there exists an existential monadic second order (EMSO) sentence such that the probability that it is true on $G(n,p=1/2)$ does not converge and conjectured that, for EMSO sentences with 2 first order variables, the 0-1 law holds. In this paper, we prove that the conjecture fails for $p\in\{\frac{3-\sqrt{5}}{2},\frac{\sqrt{5}-1}{2}\}$, and give new examples of sentences with fewer variables without convergence (even for $p=1/2$).

\vspace{0.2cm}

{\bf Keywords:} existential monadic second order language, random graphs, zero-one laws, convergence laws.

}

\vspace{0.5cm}

In this paper, we study existential monadic second order (EMSO) properties of undirected graphs. In 2001, J.-M. Le Bars proved that there exists an EMSO sentence about undirected graphs such that ${\sf P}(G_n\models\phi)$ does not converge (here, the probability distribution is uniform over the set of all graphs on the set of vertices $\{1,\ldots,n\}$). In the same paper, he conjectured that, for EMSO sentences with 2 first order variables, the 0-1 law holds (every sentence is either true asymptotically almost surely (a.a.s.), or false a.a.s.). We give an example of EMSO sentence with 1 monadic variable without convergence and an example of EMSO sentence with 3 first order variables without convergence. Moreover, we consider the binomial random graph $G(n,p)$ and move from the uniform case $p=1/2$ described above to the case of arbitrary constant $p\in(0,1)$. The above results are also true for this graph. In addition, we consider the set of EMSO sentences of the form $\exists X\,\phi(X)$ where the first order part $\phi(X)$ has quantifier depth $2$, and $X$ is the only monadic variable. We prove that, for these sentences, zero-one law holds if and only if $p\notin\{\frac{3-\sqrt{5}}{2},\frac{\sqrt{5}-1}{2}\}$. If $p\in\{\frac{3-\sqrt{5}}{2},\frac{\sqrt{5}-1}{2}\}$, there is even no convergence.

The paper is organised in the following way. Section~\ref{definition} is devoted to logical preliminaries. In Section~\ref{laws},  we review known results about logical laws of binomial random graphs. In the same section, we motivate and state new results. The remaining sections are devoted to the proofs. Their structures are described in the end of Section~\ref{laws}.

\section{Logical preliminaries}
\label{definition}

Studying zero-one laws requires an amount of logical prerequisites. We review some of the basics in this section, and refer the reader to~\cite{Logic1,ComptonLogics,Logic2,Libkin,Strange,Veresh,Survey}. Sentences in the first order language of graphs (FO sentences) are constructed using relational symbols $\sim$ (interpreted as adjacency) and $=$, logical connectives $\neg,\rightarrow,\leftrightarrow,\vee,\wedge$, variables $x,y,x_1, \ldots$ that express vertices of a graph, quantifiers $\forall,\exists$ and parentheses $($, $)$.  Monadic second order, or MSO, sentences (see~\cite{ComptonLogics,Libkin,Tobias,Zhuk_MSO}) are built of the above symbols of the first order language, as well as the variables $X,Y,X_1,\ldots$ that are interpreted as unary predicates, i.e. {\it subsets} of the vertex set. In an MSO sentence, variables $x,y,x_1, \ldots$ (that express vertices) are called {\it FO variables}, and variables $X,Y,X_1,\ldots$ (that express sets) are called {\it MSO variables}. If, in an MSO sentence $\varphi$, all the MSO variables are existential and in the beginning (that is $\varphi=\exists X_1\ldots\exists X_m\,\,\phi(X_1,\ldots,X_m)$ where $\phi(X_1,\ldots,X_m)$ is a first order sentence with unary predicates $X_1,\ldots,X_m$), then the sentence is called {\it existential} monadic second order (EMSO). Sentences must have finite length.

Following \cite{Libkin,Survey}, we call the number of nested quantifiers in a longest sequence of nested quantifiers of a sentence $\varphi$ {\it the quantifier depth} $\qd(\varphi)$ (in~\cite{Libkin}, the notion {\it quantifier rank} is used instead, but we prefer the previous one). For example, the MSO sentence
$$
 \forall X \quad\biggl[(\exists x_1\exists x_2\,\,
 [X(x_1)\wedge \neg X(x_2)])\rightarrow
 (\exists y \exists z\,\,[X(y)\wedge\neg X(z)\wedge y\sim z])\biggr]
$$
has quantifier depth $3$ and expresses the property of being connected (and its first order part has quantifier depth 2). It is known that the property of being connected {\it cannot be expressed} by a FO sentence. This fact (and many other facts about an expressibility) may be easily proved using Ehrenfeucht games (see e.g.,~\cite{Libkin,Survey}). In Section~\ref{EHR_Fagin}, we consider a modification of this game which we use in our proofs.

The quantifier depth of a sentence has the following clear algorithmic interpretation: an FO sentence of quantifier depth $k$ on an $n$-vertex graph can be verified in $O(n^k)$ time. It is very well known (see, e.g.,~\cite{Libkin}, Proposition 6.6) that the same is true for the number of variables: an FO sentence with $k$ variables on an $n$-vertex graph can be verified in $O(n^k)$ time. The later statement is stronger because, clearly, every FO sentence of quantifier depth $k$ may be rewritten using at most $k$ variables.

In what follows, for a sentence $\varphi$, we use the usual notation from model theory $G\models\varphi$ if $\varphi$ is true for $G$.

\subsection{Ehrenfeucht games}
\label{EHR_Fagin}

An important tool that allows exploiting combinatorial techniques for proving results in logic is a class of games called Ehrenfeucht games~\cite{ComptonLogics,Libkin,Strange,Veresh,Survey,Tobias,Ehren}. We consider the following modification of this game (which is also called {\it the 1,2-Fagin game}, see~\cite{Libkin}). Let $A,B$ be two graphs. There are three rounds in the game EHR$(A,B)$ and two players, called Spoiler and Duplicator. In the first round, Spoiler chooses a set of vertices $X$ in $A$, Duplicator responds with a set of vertices $Y$ in $B$. In each of the remaining rounds $i\in\{2,3\}$, Spoiler chooses either a vertex $x_{i-1}$ of $A$ or a vertex $y_{i-1}$ of $B$. Duplicator then chooses a vertex in the other graph (either $y_{i-1}$ or $x_{i-1}$). At the end of the game the vertices $x_1,x_2$ of $A$, and $y_1,y_2$ of $B$ are chosen. Duplicator wins if and only if the following properties hold

\begin{enumerate}

\item[$\cdot$] $(x_1\sim x_2)\leftrightarrow (y_1\sim y_2)$, 

\item[$\cdot$] $(x_1= x_2)\leftrightarrow (y_1=y_2)$,

\item[$\cdot$] for $i\in\{1,2\}$, $(x_i\in X)\leftrightarrow(y_i\in Y)$.

\end{enumerate}

Consider a {\it logic} (class of sentences) ${\cal E}_1^2$ containing all EMSO sentences of the form $\exists X\,\,\phi(X)$ where $\phi(X)$ is an FO sentence of quantifier depth $2$, and $X$ is the only monadic variable. The following result establishes the well-known connection between logic and the Ehrenfeucht games (we drop the proof, because it repeats usual proofs of all such results, see, e.g.,~\cite{Libkin,Veresh,Ehren,DAM}). 
\begin{theorem}
If Duplicator wins EHR$(A,B)$, then every sentence $\varphi\in\mathcal{E}_1^2$ does not differ between $A$ and $B$ (that is either $A\models\varphi$ and $B\models\varphi$, or $A\models\neg\varphi$ and $B\models\neg\varphi$).
\label{EHR}
\end{theorem}

\subsection{Extension axioms}
\label{axioms}

Another important tool that we exploit in our proofs is extension axioms.

Let $k$ be a positive integer, and $a\leq k$ be a non-negative integer. Then the sentence 
$$
 \Phi_k=\bigwedge_{a=0}^k\Phi_{k,a},\quad
 \Phi_{k,a}=\forall v_1\ldots\forall v_a\forall u_1\ldots\forall u_{k-a}
$$
$$
 \left[\left(\bigwedge_{i\in\{1,\ldots,a\},j\in\{1,\ldots,k-a\}}v_i\neq u_j\right)\rightarrow
 \left(\exists z\,\,\left[\left(\bigwedge_{i=1}^a z\sim v_i\right)\wedge\left(\bigwedge_{j=1}^{k-a}[z\nsim u_j\wedge z\neq u_j]\right)\right]\right)\right]
$$
saying that, for every set $X$ of $k$ distinct vertices and every subset $Y$ of $X$ (not necessarily non-empty), there is a vertex outside $X$ adjacent to every vertex in $Y$ and non-adjacent to every vertex in $X\setminus Y$, is called the {\it k-th extension axiom}. It is very well known (see, e.g., \cite{Libkin,Strange,Survey}) that if $A\models\Phi_k$, $B\models\Phi_k$, then any FO sentence of quantifier depth $k+1$ does not differ between $A$ and $B$.

In this paper, we also ask, is it true that, for some $k\geq 2$ and for the class $\mathcal{E}_1^k$ of EMSO sentences with only 1 monadic variable and quantifier depth $k$, there exists an $s$ such that $A\models\Phi_s$, $B\models\Phi_s$ implies that any sentence from $\mathcal{E}_1^k$ does not differ between $A$ and $B$? The negative answer (even for $k=2$) is given in Section~\ref{depth_2}.

\section{Zero-one laws and non-convergence}
\label{laws}

In 1959, P. Erd\H os and A. R\' enyi, and independently E. Gilbert, introduced two closely related models for generating random graphs. A seminal paper of Erd\H os and R\'enyi \cite{Erdos}, that appeared one year later, brought a lot of attention to the subject, giving birth to  Erd\H os-R\'enyi random graphs. In spite of the name, the more popular model $G(n,p)$  is the one proposed by Gilbert. In this model, we have $G(n,p) = (V_n,E)$, where $V_n=\{1,\ldots,n\}$, and each pair of vertices is connected by an edge with probability $p$ and independently of other pairs. For more information, we refer readers to the books \cite{Bollobas,Janson,AS}. Y. Glebskii, D. Kogan, M. Liogon'kii and V. Talanov in 1969~\cite{Glebskii}, and independently R. Fagin in 1976~\cite{Fagin}, proved that any FO sentence is either true with asymptotical probability 1 (asymptotically almost surely or a.a.s.) or a.a.s. false for $G(n,1/2)$, as $n\to \infty$. In such a situation we say that $G(n,p)$ \textit{obeys the FO zero-one law}. More generally, consider a logic $\mathcal{L}$. We say that $G(n,p)$ {\it obeys the $\mathcal{L}$ zero-one law} if, for every sentence $\varphi\in\mathcal{L}$, $\lim_{n\to\infty}{\sf P}(G(n,p)\models\varphi)\in\{0,1\}$. A weaker version of this law is called the convergence law: $G(n,p)$ {\it obeys the $\mathcal{L}$ convergence law} if, for every sentence $\varphi\in\mathcal{L}$, the limit $\lim_{n\to\infty}{\sf P}(G(n,p)\models\varphi)$ exists (but not necessarily equals 0 or 1).

To the best of our knowledge, J.F. Lynch~\cite{LynchEhrenFirst} was the first who noticed that Ehrenfeucht games may be applied for proofs of zero-one laws and J. Spencer (\cite{Spencer_Ehren}, Theorem 3.1) was the first who applied it for $G(n,p)$ for various parameters $p=p(n)$. For a logic $\mathcal{L}$, consider a version of Ehrenfeucht game (if it exists) for which an analogue of Theorem~\ref{EHR} holds (that is, if Duplicator has a winning strategy in the game on two graphs, then no sentence from $\mathcal{L}$ differs between these two graphs). Then it is not difficult to show (see, e.g.,~\cite{Survey,AS,Spencer_Ehren}) the following implication. If, for any number of rounds, a.a.s. there exists a winning strategy of Duplicator on two independent random graphs $G(n,p)$ and $G(m,p)$ (as $n,m\to\infty$), then the zero-one law for sentences from $\mathcal{L}$ holds. For FO, it is clear that (see, e.g.,~\cite{Spencer_Ehren}, Theorem 4.1) if, for an integer $k\geq 1$, $\lim_{n\to\infty}{\sf P}(G(n,p)\models\Phi_k)= 1$, then a.a.s. Duplicator wins in $k+1$ rounds. It is a simple combinatorial exercise, to prove the latter convergence (see, e.g.,~\cite{Spencer_Ehren}, Theorem 4.1). So, $G(n,p)$ obeys the FO zero-one law for every constant $p$.

However, $G(n,1/2)$ does not obey even the convergence law for EMSO. The respective construction $\varphi$ was obtained by J.-M. Le Bars in 2001~\cite{Le_Bars}. The first MSO sentence with one binary relation that has no asymptotic probability was constructed by M. Kaufmann and S. Shelah in 1985~\cite{Kaufmann_Shelah}. In 1987~\cite{Kaufmann}, Kaufmann proved that there exists an EMSO sentence with 4 binary relations that has no asymptotic probability. Note that this construction contains 4 monadic variables and 9 first order variables, and the sentence $\varphi$ proposed by Le Bars has even more variables (of both types). In the above mentioned paper, Le Bars conjectured that, for EMSO sentences with 2 first order variables, $G(n,1/2)$ obeys the zero-one law.

We disprove the conjecture for $G(n,p),$ $p\in\{\frac{3-\sqrt{5}}{2},\frac{\sqrt{5}-1}{2}\}$. More precisely, we prove the following theorem.

\begin{theorem}
Let $p\in(0,1)$ be a constant. Then $G(n,p)$ obeys the $\mathcal{E}_1^2$ zero-one law if and only if $p\notin\{\frac{3-\sqrt{5}}{2},\frac{\sqrt{5}-1}{2}\}$.

If $p\in\{\frac{3-\sqrt{5}}{2},\frac{\sqrt{5}-1}{2}\}$, then $G(n,p)$ does not obey the $\mathcal{E}_1^2$ convergence law.
\label{thm_1_2}
\end{theorem}


Moreover, we prove that 3 first order variables is enough to get non-convergence, and 1 monadic variable is also enough for this.

\begin{theorem}
Let $p\in(0,1)$ be a constant. There exist an EMSO sentence with 3 first order variables $\varphi_1$ and an EMSO sentence with 1 monadic variable $\varphi_2$, such that, for both $j=1,2$, the probability ${\sf P}(G(n,p)\models\varphi_j)$ does not converge as $n\to\infty$.
\label{thm_1_and_3}
\end{theorem}


The proofs are organized in the following way. In Section~\ref{depth_2} we state a result which answers the question from Section~\ref{axioms} and together with Lemma~\ref{main_lem} from Section~\ref{cliques_independent_sets} implies Theorem~\ref{thm_1_2}. Constructions of $\phi_1$ and $\phi_2$ from Theorem~\ref{thm_1_and_3} are given in Section~\ref{examples}.

\section{When do extension axioms imply the EMSO equivalence?}
\label{depth_2}

In Section~\ref{axioms}, we ask, is it true that, for some $k\geq 2$, there exists an $s$ such that $A\models\Phi_s$, $B\models\Phi_s$ implies that any sentence from $\mathcal{E}_1^k$ does not differ between $A$ and $B$? The answer is `no' even for $k=2$ (see Theorem~\ref{Alice} below in this section). But can we restrict somehow (in an optimal way) the set of all pairs of graphs such that the answer becomes `yes'?

We prove that there are six special monadic sentences $\varphi^1_C,\varphi^2_C,\varphi^3_C$, $\varphi^1_I,\varphi^2_I,\varphi^3_I$ (and this set is minimal, see Theorem~\ref{Alice}) such that if two graphs have the $s$-extension property (for $s$ large enough) and $\varphi^1_C,\varphi^2_C,\varphi^3_C$, $\varphi^1_I,\varphi^2_I,\varphi^3_I$ do not differ between these two graphs, then no $\varphi\in\mathcal{E}_1^2$ differs between them. These sentences are defined in the following way. Let $\max\cl(X)$ say that $X$ is a maximal clique:
$$
\max\cl(X)=[\cl(X)]\wedge[\neg\exists x\,\,(\neg X(x)\wedge[\forall y\,\,(X(y)\rightarrow x\sim y)])],
$$
where $\cl(X)=\forall x \forall y\,\,[(X(x)\wedge X(y)\wedge x\neq y)\rightarrow x\sim y]$.
Let $\varphi_C(X)$ say that there exists an external vertex which is not adjacent to any vertex of $X$:
$$
 \varphi_C(X)=\exists x\,\,(\neg X(x))\wedge(\forall y\,\,[X(y)\rightarrow x\nsim y]).
$$
Then
\begin{equation}
\varphi^1_C=\forall X\,\,[\cl(X)\rightarrow\varphi_C(X)],\quad
\varphi^2_C=\forall X\,\,[\max\cl(X)\rightarrow\neg\varphi_C(X)],
\label{phi_C}
\end{equation}
$$
\varphi^3_C=\forall X\,\,[(\cl(X)\wedge\neg\phi_C(X))\rightarrow\max\cl(X)].
$$
Note that $\varphi^1_C$ implies $\varphi^3_C$, and $\varphi^2_C$ implies $\neg\varphi^1_C$. In other words, sets 
$$
\mathcal{G}_1:=\{G:\,G\models\varphi^1_C\},
$$
$$
\mathcal{G}_2:=\{G:\,G\models\neg\varphi^1_C\wedge\neg\varphi^2_C\wedge\varphi^3_C\},
$$
$$
\mathcal{G}_3:=\{G:\,G\models\neg\varphi^1_C\wedge\neg\varphi^2_C\wedge\neg\varphi^3_C\},
$$
$$
\mathcal{G}_4:=\{G:\,G\models\varphi^2_C\wedge\varphi^3_C\},
$$
$$
\mathcal{G}_5:=\{G:\,G\models\varphi^2_C\wedge\neg\varphi^3_C\}
$$ 
form a partition of the set of all graphs.

In the same way, $\ind(X)$ say that $X$ is an independent set, $\max\ind(X)$ say that $X$ is a maximal independent set, and $\phi_I(X)$ say that there exists an external vertex which is adjacent to any vertex of $X$. Then
$$
\varphi^1_I=\forall X\,\,[\ind(X)\rightarrow\varphi_I(X)],\quad
\varphi^2_I=\forall X\,\,[\max\ind(X)\rightarrow\neg\varphi_I(X)],
$$
$$
\varphi^3_I=\forall X\,\,[(\ind(X)\wedge\neg\phi_I(X))\rightarrow\max\ind(X)].
$$
The same implications (as for $\varphi_C$-sentences) hold for these sentences.

\begin{theorem}

There exists $s$ such that if $G,H$ have the $s$-extension property, and $\varphi^1_C,\varphi^2_C,\varphi^3_C$, $\varphi^1_I,\varphi^2_I,\varphi^3_I$ do not differ between $G$ and $H$,  then no $\varphi\in\mathcal{E}_1^2$ differs between $G$ and $H$.

Moreover, the set of sentences $\varphi^1_C,\varphi^2_C,\varphi^3_C,\varphi^1_I,\varphi^2_I,\varphi^3_I$ is a minimum set such that the first part of this theorem holds in the following sense. For every positive integer $s$ and every
$$
\varphi\in\{\varphi^1_C,\neg\varphi^1_C\wedge\neg\varphi^2_C\wedge\varphi^3_C,
\neg\varphi^1_C\wedge\neg\varphi^2_C\wedge\neg\varphi^3_C,\varphi^2_C\wedge\varphi^3_C,\varphi^2_C\wedge\neg\varphi^3_C\},
$$
there exists a graph $G$ such that $G$ has the $s$-extension property and $G\models\varphi$ (similar statement holds for the sentences $\varphi^1_I,\varphi^2_I,\varphi^3_I$).

\label{Alice}
\end{theorem}

{\it Proof}. Fix $s$ as large as desired. Let $G,H$ have the $s$-extension property, and the $6$ sentences from the statement of Theorem~\ref{Alice} do not differ between $G$ and $H$. Let us prove that no sentence from $\mathcal{E}_1^2$ differs between $G$ and $H$. By Theorem~\ref{EHR}, it is sufficient to prove that Duplicator has a winning strategy in EHR$(G,H)$. In the first round, Spoiler makes a move in $G$, he chooses $X\subset V(G)$.

In what follows, we need the following definitions. Let $A$ be a graph, and $S\subset V(A)$. We say that {\it vertices $v_1,\ldots,v_s$ of $A$ have the extension property w.r.t. $S$} if for every $0\leq a\leq s$ and every pairwise distinct $i_1,\ldots,i_a\in\{1,\ldots,s\}$ there exists $z\in S$ which is adjacent to $v_{i_1},\ldots,v_{i_a}$ and non-adjacent to all the other vertices from $\{v_1,\ldots,v_s\}$. We say that {\it $A$ has the $s$-extension property w.r.t. $S$}, if any vertices $v_1,\ldots,v_s$ of $A$ have the extension property w.r.t. $S$.

1. Assume that $G$ has the $1$-extension property w.r.t. $X$ and the $1$-extension property w.r.t. $\overline{X}$.

\begin{claim}
There exist $Y\subset V(H)$ such that $H$ has the $1$-extension properties w.r.t. $Y$ and $\overline{Y}$.
\end{claim}

{\it Proof}. Let $y_1,y_2$ be vertices of $H$. Denote by $N_H(y)$ the set of all neighbors of a vertex $y$ in $H$.  Let
$$
 Y=[N_H(y_1)\cap \overline{N}_H(y_2)]\cup [\overline{N}_H(y_1)\cap N_H(y_2)].
$$
Let us prove that $H$ has the $1$-extension properties w.r.t. $Y$ and $\overline{Y}$. First, choose $y\in Y\setminus\{y_1,y_2\}$. Without loss of generality, $y\in N_H(y_1)\cap \overline{N}_H(y_2)$. By the $3$-extension property of $H$, there is a neighbor of $y$ in $Y$ (because the set $N_H(y)\cap N_H(y_1)\cap \overline{N}_H(y_2)$ is non-empty) and a non-neighbor of $y$ in $Y$ (because the set $\overline{N}_H(y)\cap N_H(y_1)\cap \overline{N}_H(y_2)\setminus\{y\}$ is non-empty). Second, choose $y\in\overline{Y}\setminus\{y_1,y_2\}$. Without loss of generality, assume that $y\in N_H(y_1)\cap N_H(y_2)$. By the $3$-extension property of $H$, there is a neighbor of $y$ in $Y$ (because the set $N_H(y)\cap N_H(y_1)\cap N_H(y_2)$ is non-empty) and there is a non-neighbor of $y$ in $\overline{Y}$ (because the set $\overline{N}_H(y)\cap N_H(y_1)\cap N_H(y_2)\setminus\{y\}$ is non-empty). Finally, let $y\in\{y_1,y_2\}$. Without loss of generality, $y=y_1$. As $N_H(y_1)\cap\overline{N}_H(y_2)\neq\varnothing$, there is a neighbor of $y_1$ in $Y$. As $\overline{N}_H(y_1)\cap N_H(y_2)\neq\varnothing$, there is a non-neighbor of $y_1$ in $Y$. As $N_H(y_1)\cap N_H(y_2)\neq\varnothing$, there is a neighbor of $y_1$ in $\overline{Y}$. As $\overline{N}_H(y_1)\cap \overline{N}_H(y_2)\setminus\{y_1\}\neq\varnothing$, there is a non-neighbor of $y_1$ in $\overline{Y}$. $\Box$\\

Let us continue describing the winning strategy of Duplicator. She chooses a set $Y\subset V(H)$ such that $H$ has the $1$-extension property w.r.t. $Y$ and the $1$-extension property w.r.t. $\overline{Y}$. Trivially, Duplicator wins in the remaining $2$ rounds.\\

2. Let $x$ be a vertex of $G$ that does not have the extension property w.r.t. $X$. Without loss of generality, assume that all neighbors of $x$ are in $\overline{X}$.

\begin{itemize}

\item[2.1] Assume that $x$ does not have the extension property w.r.t. $\overline{X}$ as well. Then $\overline{X}$ is the set of all neighbors of $x$ (or $x$ is an extra vertex of $\overline{X}$). Obviously, in this case, every other vertex of $G$ has the extension property w.r.t. $X$ and w.r.t. $\overline{X}$. This follows from the 2-extension property of $G$: for every $\tilde x\neq x$, the sets $N_G(x)\cap N_G(\tilde x)$, $N_G(x)\cap \overline{N}_G(\tilde x)\setminus\{\tilde x\}$, $\overline{N}_G(x)\cap N_G(\tilde x)\setminus\{x\}$, $\overline{N}_G(x)\cap \overline{N}_G(\tilde x)\setminus\{x,\tilde x\}$ are not empty. Let $y$ be a vertex of $H$. Duplicator chooses $Y=\overline{N}_H(y)$ (or $Y=\overline{N}_H(y)\setminus\{y\}$, whenever $x\in\overline{X}$). By the $2$-extension property, every $\tilde y\neq y$ has the extension property w.r.t. $Y$ and w.r.t. $\overline{Y}$. Trivially, Duplicator wins in the remaining $2$ rounds.

\item[2.2] Let $x$ have the extension property w.r.t. $\overline{X}$. As above, all the other vertices of $G$ have the extension property w.r.t. $\overline{X}$.

    Assume that all the vertices of $X$ have the extension property w.r.t. $X$. In this case, $x\in\overline{X}$. If there are no vertices which are adjacent to all vertices in $X$, take arbitrary vertex $y$ in $H$ and set $Y=\overline{N}_H(y)\setminus\{y,\tilde y\}$, where $\tilde y$ is an arbitrary non-neighbor of $y$. If there is a vertex $\hat y$ which is adjacent to all vertices in $Y$, then $\overline{N}_H(y)\cap N_H(\tilde y)\cap\overline{N}_H(\hat y)=\varnothing$. This contradicts the $3$-extension property of $H$. If, in $Y$, there is a vertex which is non-adjacent to all the other vertices in $Y$, then $\overline{N}_H(y)\cap N_H(\tilde y)\cap N_H(\hat y)=\varnothing$. Again, this contradicts the 3-extension property. Duplicator chooses $Y$. If an $\tilde x$ (in $\overline{X}$) is a common neighbor of all the vertices in $X$, then Duplicator chooses $Y$ such that $H|_Y=P_2\sqcup P_2$. By the $4$-extension property, such an $Y$ exists, and there is a common neighbor and a common non-neighbor of all the vertices of $Y$ in $H$.

    Let $G|_X$ be a clique. In this case, $x\in\overline{X}$ as well. If there are no vertices which are adjacent to all the vertices in $X$, then $G\models\neg\phi^2_C$. As $\phi^2_C$ does not differ between $G$ and $H$, $\phi^2_C$ is also false on $H$. Therefore, there exists $Y\subset V(H)$ such that $H|_Y$ is a maximal clique, and there exists a vertex $y\in\overline{Y}$ which is not adjacent to any vertex of $Y$. It is easy to see that every vertex of $H$ has the extension property w.r.t. $\overline{Y}$ (otherwise, we easily get a contradiction with the 2-extension property). Duplicator chooses $Y$. If $\tilde x$ is a common neighbor of all vertices in $X$, then Duplicator chooses $Y=\{y_1,y_2\}$ such that $y_1\sim y_2$. By the $3$-extension property, such an $Y$ exists, and there exist vertices $y,\tilde y\in\overline{Y}$ such that $y$ is not adjacent to any of $y_1,y_2$, and $\tilde y$ is a common neighbor of $y_1,y_2$. It is easy to see that every vertex of $H$ has the extension property w.r.t. $\overline{Y}$.

    If $X$ is an independent set, then one can consider four cases: in $\overline{X}$, 1) there is a common neighbor of all vertices of $X$ and a vertex which is not adjacent to any vertex of $X$, 2) there is a common neighbor of all vertices of $X$ and there are no vertices which are not adjacent to any vertex of $X$, 3) there are no common neighbors of all vertices of $X$ and there is a vertex which is not adjacent to any vertex of $X$, 4) there are no common neighbors of all vertices of $X$ and no vertices which are not adjacent to any vertex of $X$. In the first case, the target set $Y$ exists due to the $3$-extension property of $H$. In the second case, the target set $Y$ exists because $H\models(\neg\phi^2_C)$. In the third case, the target set $Y$ exists because $H\models(\neg\phi^1_C)\wedge(\neg\phi^3_C)$. In the fourth case, the target set $Y$ exists because $H\models(\neg\phi^1_C)$.

    If $G|_X$ is not a clique but in $X$ there is a vertex $\tilde x$ such that $[N_G(\tilde x)\cup \{\tilde x\}]\supset X$, then $x\in\overline{X}$. If there is also a vertex $\hat x\in\overline{X}$ which is a common neighbor of all vertices in $X$, then, by the $4$-extension property, there is a set $Y=\{\tilde y,y_1,y_2\}$ and vertices $y,\hat y$ in $H$ such that $y$ is a common neighbor of $y_1,y_2$, $\hat y$ is a common neighbor of $\tilde y,y_1,y_2$, $y$ is non-adjacent to any vertex of $Y$, and $y_1\nsim y_2$. Duplicator chooses the set $Y$. If, in $\overline{X}$, there are no common neighbors of all vertices in $X$, then find two vertices $\tilde y$ and $y$ that are not adjacent in $H$. Let $Y=[N_H(\tilde y)\cap\overline{N}_H(y)]\cup\{\tilde y\}$. Suppose that there is a vertex $\hat y\in\overline{Y}$ which is a common neighbor of all vertices in $Y$. Then $N_H(\tilde y)\cap\overline{N}_H(y)\cap \overline{N}_H(\hat y)=\varnothing$. This contradicts the $3$-extension property of $H$.

    Finally, let $X$ be not an independent set and there be an isolated vertex in $G|_X$. Consider four cases: in $\overline{X}$, 1) there is a common neighbor of all vertices of $X$ and a vertex which is not adjacent to any vertex of $X$, 2) there is a common neighbor of all vertices of $X$ and there are no vertices which are not adjacent to any vertex of $X$, 3) there are no common neighbors of all vertices of $X$ and there is a vertex which is not adjacent to any vertex of $X$, 4) there are no common neighbors of all vertices of $X$ and no vertices which are not adjacent to any vertex of $X$. In the first case, the target set $Y$ exists due to the $4$-extension property of $H$. In the second case, find two vertices $\tilde y$ and $y$ that are adjacent in $H$. Let $Y=\overline{N}_H(y)\cap N_H(\tilde y)$. Suppose that there is a vertex $\hat y\in\overline{Y}$ which is not adjacent to any vertex in $Y$. Then $(\overline{N}_H(y)\cap N_H(\tilde y)\cap N_H(\hat y))\setminus\{y\}=\varnothing$. This contradicts with the $3$-extension property of $Y$. The existence of the target set $Y$ in the third case follows from the $3$-extension property in same way. In the fourth case, find two vertices $\tilde y$ and $y$ that are non-adjacent in $H$. Set $Y=\overline{N}_H(y)\setminus\{\tilde y\}$. It easily follows from the $3$-extension property that there are no common neighbors of all vertices of $Y$ and no vertices which are not adjacent to any vertex of $Y$.

    In all the above cases, Duplicator has a winning strategy in the remaining two rounds.

\end{itemize}

Now, let us prove the second part of the theorem. Obviously, it is enough two prove it only for sentences $\varphi^1_C,\varphi^2_C,\varphi^3_C$. For $\varphi\in\{\varphi^1_C,\neg\varphi^1_C\wedge\neg\varphi^2_C\wedge\neg\varphi^3_C\}$, the existence of the required graph follows from Lemma~\ref{main_lem} of Section~\ref{cliques_independent_sets}.

Let $\varphi=\varphi^2_C\wedge\varphi^3_C$. Consider $\ell$ complete bipartite graphs $H_i\cong K_{n,n}$ and draw an edge with probability $p$ between every two vertices from different graphs (all edges appear independently). Denote this random graph on $2n\ell$ vertices by $G_n^1$. Set $V(H_i)=V_i$, $i\in\{1,\ldots,\ell\}$. It can be easily proved that a.a.s. (as $n\to\infty$) this graph has the $(\ell-1)$-extension property. Moreover, a.a.s. every maximal clique has size $2\ell$ and has two common vertices with every $H_i$. For every clique $S$ with $|V(S)|\leq 2\ell-1$, there is $V_i$ such that $|V_i\cap V(S)|\leq 1$. In $V_i$, there are at least $n-1$ vertices that are not adjacent to the only vertex (if it exists) in $V_i\cap V(S)$. Therefore, a.a.s. for every such a clique $S$, there is a vertex which has no neighbors in $S$. But for every clique $S$ of size $2\ell$ and every $i$, each vertex in $V_i\setminus V(S)$ has a neighbor in $V_i\cap V(S)$. So, a.a.s. $G_n^1\models\varphi$.

Let $\varphi=\neg\varphi^1_C\wedge\neg\varphi^2_C\wedge\varphi^3_C$. Consider $\ell-1$ complete bipartite graphs $H_i\cong K_{n,n}$ and one star graph $H_{\ell}=K_{1,n}$. Let $x$ be the central vertex of $H_{\ell}$. Draw an edge with probability $p$ between every two vertices from different graphs (all edges appear independently). Denote this random graph on $2n(\ell-1)+n+1$ vertices by $G^2_n$. Set $V(H_i)=V_i$, $i\in\{1,\ldots,\ell\}$. As above, a.a.s. this graph has the $(\ell-1)$-extension property. Moreover, a.a.s. every maximal clique has either size $2\ell-1$ and does not contain $x$, or has size $2\ell$ and contains $x$ (both types of cliques appear in the random graph a.a.s.). In $G_n^2$, for every clique $S$ of size $2\ell$ and every vertex $u\notin V(S)$, $u$ has a neighbor in $V(S)$. Moreover, a.a.s. every clique $S$ of size $2\ell-1$ has an outside vertex which has no neighbors in $V(S)$. Therefore, a.a.s. $G_n^2\models\varphi$.

Finally, let $\varphi=\varphi^2_C\wedge\neg\varphi^3_C$. Consider $\ell-1$ complete bipartite graphs $H_i\cong K_{n,n}$ and one union $H_{\ell}$ of a complete bipartite graph $K_{n,n}$ with a vertex which is adjacent to all $2n$ vertices of $K_{n,n}$. Let $x$ be the universal vertex of $H_{\ell}$. Draw an edge with probability $p$ between every two vertices from different graphs (all edges appear independently), and denote this random graph by $G^3_n$. Set $V(H_i)=V_i$, $i\in\{1,\ldots,\ell\}$. As above, a.a.s. $G^3_n$ has the $(\ell-1)$-extension property. Moreover, a.a.s. every maximal clique has either size $2\ell$ and does not contain $x$, or has size $2\ell+1$ and contains $x$ (both types of cliques appear in the random graph a.a.s.). Note that, for every $2\ell$-subgraph $S$ of a maximal $2\ell+1$-clique of $G^3_n$ which has two common vertices with each $V_i$, each $u\notin S$ has a neighbor in $S$. Therefore, a.a.s. $G_n^3\models\varphi$. $\Box$.

\section{Maximal independent sets and maximal cliques in the random graph}
\label{cliques_independent_sets}

In this section, we ask, which of the sets $\mathcal{G}_j$, $j\in\{1,2,3,4,5\}$, the random graph $G(n,p)$ belongs a.a.s. for different values of $p$. Surprisingly, for $p=\frac{3-\sqrt{5}}{2}$, there is no such a set. Moreover, for this value of $p$, the probability ${\sf P}(G(n,p)\models\varphi^3_C)$ does not converge. A similar result holds true for $\varphi_I$ sentences.

\begin{lemma}
Let $p=\const\in(0,1)$.

\begin{enumerate}

\item If $p>\frac{3-\sqrt{5}}{2}$, then ${\sf P}(G(n,p)\models(\neg\varphi^1_C)\wedge(\neg\varphi^2_C)\wedge(\neg\varphi^3_C))\to 1$ as $n\to\infty$.

\item If $p<\frac{3-\sqrt{5}}{2}$, then ${\sf P}(G(n,p)\models\varphi^1_C)\to 1$ as $n\to\infty$.

\item If $p=\frac{3-\sqrt{5}}{2}$, then ${\sf P}(G(n,p)\models\varphi^3_C)$ does not converge.

\item If $p<\frac{\sqrt{5}-1}{2}$, then ${\sf P}(G(n,p)\models(\neg\varphi^1_I)\wedge(\neg\varphi^2_I)\wedge(\neg\varphi^3_I))\to 1$ as $n\to\infty$.

\item If $p>\frac{\sqrt{5}-1}{2}$, then ${\sf P}(G(n,p)\models\varphi^1_I)\to 1$ as $n\to\infty$.

\item If $p=\frac{\sqrt{5}-1}{2}$, then ${\sf P}(G(n,p)\models\varphi^3_I)$ does not converge.

\end{enumerate}
\label{main_lem}
\end{lemma}

{\it Proof}. It is enough to prove 1, 2 and 3. Let $k$ be a positive integer, and $X_1(k)$ be the number of pairs $(C_1,c_1)$ where $C_1$ is a maximal clique of size $k$ in $G(n,p)$ and $c_1$ is a vertex which is not adjacent to any vertex of $C_1$. Let $X_2(k)$ be the number of pairs $(C_2,c_2)$ where $C_2$ is a clique of size $k$ in $G(n,p)$, $c_2$ is a (external) common neighbor of all vertices of $C_2$, and there does not exist a vertex which is not adjacent to any vertex of $C_2$. We use Chebyshev's inequality to prove the lemma. First,
$$
 {\sf E}X_1(k)={n\choose k}(n-k)p^{{k\choose 2}}(1-p)^k(1-p^k)^{n-k-1},
$$
\begin{equation}
 {\sf E}X_2(k)={n\choose k}(n-k)p^{k\choose 2}p^k(1-(1-p)^k)^{n-k-1}.
\label{expectation_x2}
\end{equation}
Second,
$$
 {\sf D}X_1(k)\leq{\sf E}X_k^1-({\sf E}X_k^1)^2+
$$
$$
 {n\choose k}(n-k)(n-k-1)p^{k\choose 2}(1-p)^{2k}+{n\choose k}{n-k\choose k}p^{2{k\choose 2}}(n-2k)^2(1-p)^{2k}(1-p^k)^{2(n-2k-2)}
$$
$$
 \sum_{\ell=1}^{k-1}{n\choose k}{k\choose \ell}{n-k\choose k-\ell}p^{2{k\choose 2}-{\ell\choose 2}}(n-2k+\ell)(1-p)^{2k}[(1-p)^{-\ell}+(n-2k+\ell-1)].
$$
$$
 {\sf D}X_2(k)\leq{\sf E}X_k^2-({\sf E}X_k^2)^2+
$$
$$
 {n\choose k}(n-k)(n-k-1)p^{k\choose 2}p^{2k}+{n\choose k}{n-k\choose k}p^{2{k\choose 2}} p^{2k}(n^2-3kn+3k^2) (1-(1-p)^k)^{2(n-k-2)}+
$$
$$
 \sum_{\ell=1}^{k-1}{n\choose k}{k\choose \ell}{n-k\choose k-\ell}p^{2{k\choose 2}-{\ell\choose 2}}\times
$$ 
$$
 \left([n-2k+\ell]p^{2k}[(k-\ell+1)p^{-\ell}+(n-2k+\ell-1)]+[k-\ell]^2p^{2k-2\ell}\right).
$$

1) Let $p\geq 1/2$. Set $k=\lfloor\log_{1/p}n\rfloor$. Then  ${\sf E}X_i(k)=e^{[\ln^2 n/(2\ln (1/p))](1+o(1))}\to\infty$ as $n\to\infty$ for both $i=1,2$.

For $i=1$, as $(1-p^k)^{n}=e^{-1}+o(1)$, we get
$$
 {\sf D}X_1(k)\leq
 {\sf E}X_1(k)+\left[{n\choose k}(n-k)(n-k-1)p^{k\choose 2}(1-p)^{2k}\right]+
 \left[{n\choose k}p^{2{k\choose 2}}(1-p)^{2k}\sum_{\ell=1}^{k-1}F_{\ell}\right]+
$$
$$
 {n\choose k}{n-k\choose k}p^{2{k\choose 2}}(n-2k)^2(1-p)^{2k}(1-p^k)^{2(n-2k-2)}
 -({\sf E}X_k^1)^2=
$$
$$
 {n\choose k}p^{2{k\choose 2}}(1-p)^{2k}\sum_{\ell=1}^{k-1}F_{\ell}+
 o\left(({\sf E}X_k^1)^2\right).
$$
where
$F_{\ell}=A_{\ell}(n-2k+\ell)[(1-p)^{-\ell}+(n-2k+\ell-1)]$,
\begin{equation}
A_{\ell}={k\choose \ell}{n-k\choose k-\ell}p^{-{\ell\choose 2}}.
\label{A_ell}
\end{equation}
By Chebyshev's inequality,
$$
 {\sf P}(X_1(k)=0)\leq\frac{{\sf D}X_1(k)}{({\sf E}X_1(k))^2}.
$$
So, it is enough to prove that $\sum_{\ell=1}^{k-1}F_{\ell}=o\left({n\choose k}(n-k)^2\right)$. The last equality follows from the fact that $F_{\ell}$ first decreases and after that increases on $\{1,\ldots,k-1\}$ (for $n$ large enough) and
$$
 kF_1=k^2{n-k\choose k-1}(n-2k+1)[(1-p)^{-1}+(n-2k)]=nk^3{n-k\choose k}(1+o(1))=o\left({n\choose k}(n-k)^2\right),
$$
$$
 k F_{k-1}=k^2(n-k)p^{-{k-1\choose 2}}(n-k-1)((1-p)^{1-k}+n-k-2)=e^{[\ln^2 n/(2\ln(1/p))(1+o(1))]}=o\left({n\choose k}\right).
$$

For $i=2$, similarly, we get
$$
 {\sf D}X_2(k)\leq
 {\sf E}X_2(k)+{n\choose k}(n-k)(n-k-1)p^{k\choose 2}p^{2k}+
 \left[{n\choose k}p^{2{k\choose 2}}p^{2k}\sum_{\ell=1}^{k-1}G_{\ell}\right]+
$$
$$
 {n\choose k}{n-k\choose k}p^{2k}(n^2-3kn+3k^2)(1-(1-p)^k)^{2(n-2k-2)}
 -({\sf E}X_2(k))^2=
$$
$$
 {n\choose k}p^{2{k\choose 2}}p^{2k}\sum_{\ell=1}^{k-1}G_{\ell}+
 o\left(({\sf E}X_2(k))^2\right),
$$
where
$G_{\ell}=A_{\ell}\left([n-2k+\ell][(k-\ell+1)p^{-\ell}+(n-2k+\ell-1)]+[k-\ell]^2p^{-2\ell}\right)$. Obviously, $G_{\ell}\leq A_{\ell}n^2k^2=:\tilde G_{\ell}$ for all $\ell\in\{1,\ldots,k-1\}$. The equality ${\sf P}(X_2(k)=0)=o(1)$ holds because $A_{\ell}$ (and, therefore, $\tilde G_{\ell}$) first decreases and then increases on $\{1,\ldots,k-1\}$ (for $n$ large enough) and
$$
 k\tilde G_1\cdot{n\choose k}p^{2{k\choose 2}}p^{2k}=
$$
$$ 
 k^4{n-k\choose k-1}{n\choose k}p^{2{k\choose 2}}p^{2k}n^2=
 \frac{k^5n^2}{n-2k+1}{n-k\choose k}{n\choose k}p^{2{k\choose 2}}p^{2k}
 =o\left(({\sf E}X_2(k))^2\right),
$$
$$
 \frac{k\tilde G_{k-1}{n\choose k}p^{2{k\choose 2}}p^{2k}}{({\sf E}X_2(k))^2}=k^4 n\frac{p^{-{{{k-1}\choose 2}}}}{{n\choose k}}[1+o(1)]=e^{-\frac{\ln^2 n}{\ln(1/p)}(1+o(1))}=o(1).
$$

Below, we consider $p<1/2$. The crucial thing is a sign of the expression $2\ln\frac{1}{1-p}-\ln\frac{1}{p}$. We will show that the expectations of the considered random variables (in each of the below three cases) equal $\exp\left[c(p)\ln^2 n\left(2\ln\frac{1}{1-p}-\ln\frac{1}{p}\right)(1+o(1))\right]$ for some constants $c(p)$. Then, our intuition is the following: if the expectation tends to infinity, then a.a.s. the respective structure exists. If it tends to 0, then a.a.s. it does not exist. And the only situation when the non-covergence is possible is when $2\ln\frac{1}{1-p}-\ln\frac{1}{p}=0$, i.e. $p=\frac{3-\sqrt{5}}{2}$.\\ 

2) Let $\frac{3-\sqrt{5}}{2}<p<1/2$. In this case, $2\ln\frac{1}{1-p}-\ln\frac{1}{p}>0$. Set $k=\lfloor\log_{1/(1-p)}n\rfloor$. Then
$$
{\sf E}X_i(k)=e^{\ln^2 n\frac{2\ln(1/(1-p))-\ln(1/p)}{2\ln^2(1/(1-p))}(1+o(1))}\to\infty\text{ as }n\to\infty
$$
for both $i=1,2$.

As above, $A_{\ell}$ (defined in~(\ref{A_ell})) first decreases and then increases on $\{1,\ldots,k-1\}$ (for $n$ large enough). Moreover, for $n$ large enough, $A_1>A_{k-1}$. Indeed,
$$
 A_1=k{n-k\choose k-1}=e^{\ln^2 n\frac{1}{\ln(1/(1-p))}(1+o(1))}>
 e^{\ln^2 n\frac{\ln(1/p)}{2\ln^2(1/(1-p))}(1+o(1))},
$$
$$
 A_{k-1}=k(n-k)p^{-{k-1\choose 2}}=e^{\ln^2 n\frac{\ln(1/p)}{2\ln^2(1/(1-p))}(1+o(1))}.
$$

For $i=1$, as $(1-p^k)^{n}=1+o(1)$, we get
$$
 {\sf D}X_1(k)\leq {n\choose k}p^{2{k\choose 2}}(1-p)^{2k}\sum_{\ell=1}^{k-1}F_{\ell}+
 o\left(({\sf E}X_1(k))^2\right).
$$
Note that, for all $\ell\in\{1,\ldots,k-1\}$, $F_{\ell}\leq 2n^2 A_{\ell}=:\tilde F_{\ell}$. As, for $n$ large enough, $\tilde F_1>\tilde F_{k-1}$, by Chebyshev's inequality, we get
$$
 {\sf P}(X_1(k)=0)\leq\frac{{\sf D}X_1(k)}{({\sf E}X_1(k))^2}\leq
 \frac{2n^2k^2{n-k\choose k-1}{n\choose k}p^{2{k\choose 2}}(1-p)^{2k}}{({\sf E}X_1(k))^2}+o(1)=o(1).
$$

For $i=2$, as $(1-(1-p)^k)^{n-k}=e^{-1}+o(1)$, similarly, we get
$$
 {\sf D}X_2(k)\leq
 {n\choose k}p^{2{k\choose 2}}p^{2k}\sum_{\ell=1}^{k-1}\tilde G_{\ell}+
 o\left(({\sf E}X_2(k))^2\right).
$$
As, for $n$ large enough, $\tilde G_1>\tilde G_{k-1}$, by Chebyshev's inequality, we get that
$$
{\sf P}(X_2(k)=0)\leq\frac{k \tilde G_1}{{n\choose k}(n-k)^2(e^{-2}+o(1))}+o(1)=O\left(\frac{k^4n^2{n-k\choose k-1}}{{n\choose k}(n-k)^2}\right)+o(1)=o(1).
$$

3) Let $p<\frac{3-\sqrt{5}}{2}$. In this case, $2\ln\frac{1}{1-p}-\ln\frac{1}{p}<0$. Let $\tilde X_2(k)$ be the number of maximal cliques $C$ such that there are no vertices which are not adjacent to any vertex of $C$. Here, we prove that $\sum_k{\sf E}\tilde X_2(k)\to 0$. From Markov inequality, this implies that ${\sf P}(\forall k\,\,\tilde X_2(k)=0)\to 1$.
We get
$$
 {\sf E}\tilde X_2(k)\leq e^{k\ln n-k\ln k+k-\frac{k^2}{2}\ln(1/p)+\frac{k}{2}\ln(1/p)
 -np^k-n(1-p)^k+o(1)}.
$$
Let us find the maximum of the function $f(k)=k\ln n-k\ln k+k-\frac{k^2}{2}\ln(1/p)+\frac{k}{2}\ln(1/p)
 -np^k-n(1-p)^k$. The function $f^{\prime}(k)=\ln n+\frac{1}{2}\ln(1/p)-\ln k-k\ln(1/p)+np^k\ln(1/p)+n(1-p)^k\ln(1/(1-p))$ decreases and equals $0$ for $k=\log_{1/(1-p)}n-\log_{1/(1-p)}\ln n(1+o(1))$. For such $k$,
$\ln{\sf E}\tilde X_2(k)\leq \frac{2\ln\frac{1}{1-p}-\ln\frac{1}{p}}{2\ln^2\frac{1}{1-p}}\ln^2 n(1+o(1))$. Finally, we get
$$
\sum_k{\sf E}\tilde X_2(k)\leq n e^{\frac{2\ln\frac{1}{1-p}-\ln\frac{1}{p}}{2\ln^2\frac{1}{1-p}}\ln^2 n(1+o(1))}=o(1).
$$

4) Let $p=\frac{3-\sqrt{5}}{2}$. First, let us estimate ${\sf E}X_2(k)$ from above. From~(\ref{expectation_x2}), we get
$$
 {\sf E}X_2(k)\leq e^{k\ln n-k\ln k+k+\ln n-k^2\ln(1/(1-p))-k\ln(1/(1-p))-n(1-p)^k}.
$$
Let us find a maximum of the function $f(k)=k\ln n-k\ln k+k-k^2\ln(1/(1-p))-k\ln(1/(1-p))-n(1-p)^k$ on $[0,n]$. The function $f^{\prime}(k)=\ln n-\ln k-2k\ln(1/(1-p))-\ln(1/(1-p))+n(1-p)^k\ln(1/(1-p))$ decreases and equals $0$ for
\begin{equation}
k^*=\frac{1}{\ln(1/(1-p))}\left(\ln n-\ln\ln n+\ln\ln\frac{1}{1-p}+\frac{\ln\ln n}{\ln n}+O\left(\frac{1}{\ln n}\right)\right).
\label{k-star}
\end{equation}
As $f(k^*)=\ln\ln n+O(1)$, we get
$$
 {\sf E}X_2(k^*)=\frac{1}{\sqrt{2\pi k}}e^{f(k^*)+o(1)}=e^{\frac{1}{2}\ln\ln n+O(1)}
$$
if $k^*$ is an integer. 

Now, let $\tilde k$ be equal to an integer in $\{1,\ldots,n\}$ for which the value of the function $f$ is maximal. We will construct two sequences of integers $n$ such that, for one of them, there exist integers $k$ satisfying~(\ref{k-star}), and, for the second one, elements of the sequence defined in~(\ref{k-star}) are close to $k+1/2$ for some integers $k$.

Note that, for $n=\lfloor k(1/(1-p))^k\rfloor$, we easily get that $k$ satisfies the equality~(\ref{k-star}). Then, for such $n$, $\tilde k$ satisfies the equality~(\ref{k-star}) as well, and so, ${\sf E}X_2(\tilde k)=e^{\frac{1}{2}\ln\ln n+O(1)}\to\infty$ as $n\to\infty.$

But, for $n=\lfloor k(1/(1-p))^{k+1/2}\rfloor$ and every function $k^*$ defined in~(\ref{k-star}), $|k-k^*|=\frac{1}{2}+o(1)$. So, for $n$ large enough, $|\tilde k-k^*|>\frac{1}{3}$. Moreover, for any $\varepsilon>0$,
$$
 f(k^*+\varepsilon)=f(k^*)-\ln n\left(\varepsilon-\frac{1-(1-p)^{\varepsilon}}{\ln(1/(1-p))}\right),\quad f(k^*-\varepsilon)=f(k^*)-\ln n\left(\frac{1/(1-p)^{\varepsilon}-1}{\ln(1/(1-p))}-\varepsilon\right).
$$
It means that, for $n$ large enough, there is a constant $c>0$ such that $f(\tilde k)\leq e^{-c\ln n}$. Finally, for some constant $C>0$,
$$
 \max_{|k-k^*|\geq C}|f(k)|\leq -2\ln n,
$$
and so, for the considered large enough $n$, 
$$
{\sf P}(\exists k\quad X_2(k)>0)\leq\sum_{k=1}^n{\sf E}X_2(k)\leq ne^{-2\ln n}+O\left(e^{-f(\tilde k)}\right)=o(1).
$$

To disprove the convergence, it remains to prove that, for the first sequence $n=\lfloor k(1/(1-p))^k\rfloor$, $\frac{{\sf D}X_2(k)}{({\sf E}X_2(k))^2}=o(1)$. Let us compute the variance. For two $k$-sets $S_1,S_2\subset\{1,\ldots,n\}$ such that $|S_1\cap S_2|=\ell$, denote $\beta_{\ell}$ the probability that a fixed vertex in $\{1,\ldots,n\}\setminus(S_1\cup S_2)$ has a neighbor in $S_1\cup S_2$. Then
\begin{equation}
 {\sf D}X_2(k)={\sf E}X_2(k)+\left[{n\choose k}(n-k)(n-k-1)p^{{k\choose 2}}p^{2k}\beta_k^{n-k-2}\right]+
\label{non_conv_variance_1}
\end{equation}
\begin{equation} 
 \left[{n\choose k}{n-k\choose k}(n-2k)^2p^{2{k\choose 2}}p^{2k}\beta_0^{n-2k-2}-({\sf E}X_2(k))^2\right]+
\label{non_conv_variance_2}
\end{equation}
\begin{equation}
 +\left[{n\choose k}p^{2{k\choose 2}}\sum_{k=1}^{\ell-1}(B^1_{\ell}+B^2_{\ell}+B^3_{\ell}+B^4_{\ell})\right],
\label{non_conv_variance_3}
\end{equation}
where the factor $B^1_{\ell}:=A_{\ell}(k-\ell)^2p^{2(k-\ell)}\beta_{\ell}^{n-2k+\ell}$ corresponds to the case when each one of two $k$-sets has the common neighbor inside the second one; the factor $B^2_{\ell}:=A_{\ell}2(k-\ell)(n-2k+\ell)p^{2k-\ell}\beta_{\ell}^{n-2k+\ell-1}$ corresponds to the case when exactly one of two $k$-sets has the common neighbor inside the second one; the factor $B^3_{\ell}:=A_{\ell}(n-2k+\ell)p^{2k-\ell}\beta_{\ell}^{n-2k+\ell-1}$ corresponds to the case when both $k$-sets have the same common; and the factor $B^4_{\ell}:=A_{\ell}(n-2k+\ell)(n-2k+\ell-1)p^{2k}\beta_{\ell}^{n-2k+\ell-2}$ corresponds to the case when $k$-sets have distinct common neighbors outside their union.
 
Obviously, for every $\ell\in\{0,1,\ldots,k-1\}$, $\beta_{\ell}=1-2(1-p)^{k}+(1-p)^{2k-\ell}$. Therefore,

\begin{center}
 for every $\ell=\mathrm{const}$ and $\phi(k)=O(k)$,
 
 $\beta_{\ell}^{n-\phi(k)}\sim \beta_k^{2n},\quad$
 $\beta_{k-\ell}^{n-\phi(k)}\sim [1-(1-p)^k(2-(1-p)^{\ell})]^n$.
\end{center}

So, the first summand in~(\ref{non_conv_variance_1})--(\ref{non_conv_variance_3}) equals ${\sf E}X_2(k)=o([{\sf E}X_2(k)]^2)$, the second summand equals $np^k{\sf E}X_2(k)(1+o(1))\sim e^{-\ln n+\frac{3}{2}\ln\ln n+o(1)}=o(1)$,  and the third summand equals $o([{\sf E}X_2(k)]^2)$ as well. So, it remains to prove that 
$$
\sum_{k=1}^{\ell-1}\max\{B_{\ell}^1,B_{\ell}^2,B_{\ell}^3,B_{\ell}^4\}=o\left({n\choose k}n^2p^{2k}\beta_k^{2n}\right).
$$
Let $G(n)={n\choose k}n^2p^{2k}\beta_k^{2n}$. As $B_{\ell}^2\geq B_{\ell}^3$ for all $\ell\in\{1,\ldots,k-1\}$, it remains to estimate from above $B_{\ell}^1,B_{\ell}^2,B_{\ell}^4$ only.

If $\ell=\mathrm{const}$, then
$$
 B_{\ell}^1\sim \frac{k^{2\ell+2}}{n^{\ell}}{n\choose k}p^{2k}\beta_k^{2n},\quad
 B_{\ell}^2\sim \frac{2k^{2\ell+1}}{n^{\ell-1}}{n\choose k}p^{2k}\beta_k^{2n},\quad
 B_{\ell}^4\sim \frac{k^{2\ell}}{n^{\ell-2}}{n\choose k}p^{2k}\beta_k^{2n};
$$
$$
 B_{k-\ell}^1=\frac{k^{2\ell}}{n^{\ell+o(1)}}p^{-k^2/2}p^{\ell k-k/2}e^{-n(1-p)^k(2-(1-p)^{\ell})}={n\choose k}p^{\ell k}k^{2\ell}n^{-\ell+1-(1-(1-p)^{\ell})/\ln[1/(1-p)]+o(1)}\leq 
$$
$$ 
 G(n)\cdot n^{-p/\ln[1/(1-p)]+o(1)},
$$
$$ 
 B_{k-\ell}^2=B_{k-\ell}^1n^{1+o(1)}p^k\leq G(n)\cdot n^{-1-p/\ln[1/(1-p)]+o(1)},
$$ 
$$ 
 B_{k-\ell}^4=B_{k-\ell}^1n^{2+o(1)}p^{2k}\leq G(n)\cdot n^{-2-p/\ln[1/(1-p)]+o(1)}.
$$
So, for $n$ large enough, 
$$
\max\{B_1^1,B_1^2,B_1^4\}\sim G(n)\cdot\frac{k^{2\ell}}{n}.
$$
Moreover,
$$
\max\{B_{k-1}^1,B_{k-1}^2,B_{k-1}^4\}\leq G(n)\cdot n^{-p/\ln[1/(1-p)]+o(1)};\quad \max\{B_{k-2}^1,B_{k-2}^2,B_{k-2}^4\}\leq G(n)\cdot n^{-3+o(1)}.
$$

It is easy to see that $A_{\ell}(k-\ell)^2p^{-2\ell}$, $A_{\ell}(k-\ell)p^{-\ell}$ and $A_{\ell}$ first decrease and then increase on $\{1,\ldots,k-1\}$. Moreover, for all $\ell\in\{1,\ldots,k-1\}$, $\beta_{\ell}\leq\beta_k$. Therefore, for all $\ell\in\{1,\ldots,k-2\}$,
$$
 \max\{B_{\ell}^1,B_{\ell}^2,B_{\ell}^4\}=O\left(\max\left\{\max\{B_1^1,B_1^2,B_1^4\}\frac{\beta_k^n}{\beta_1^n},
 \max\{B_{k-2}^1,B_{k-2}^2,B_{k-2}^4\}\frac{\beta_k^n}{\beta_{k-2}^n}\right\}\right)=
$$ 
$$ 
 =O\left(G(n)\cdot\max\left\{n^{-(1-1/\ln[1/(1-p)])+o(1)},
 n^{-(3-(2p-p^2)/\ln[1/(1-p)])+o(1)}\right\}\right)=n^{-(1-1/\ln[1/(1-p)])+o(1)}.
$$
Finally, we get
$$
 \sum_{k=1}^{\ell-1}\max\{B_{\ell}^1,B_{\ell}^2,B_{\ell}^3,B_{\ell}^4\}=
 G(n)\left[n^{-p/\ln[1/(1-p)]+o(1)}+kn^{-(1-1/\ln[1/(1-p)])+o(1)}\right]=
 o\left(G(n)\right).\quad\Box
$$


Let $p=\const\in(0,1)\setminus\{\frac{3-\sqrt{5}}{2},\frac{\sqrt{5}-1}{2}\}$. From Theorem~\ref{Alice} and Lemma~\ref{main_lem}, it easily follows that, for every $\mathcal{E}_1^2$ sentence $\varphi$, the limit $\lim_{n\to\infty}{\sf P}(G(n,p)\models\varphi)$ exists and equals either $0$ or $1$. If $p\in\{\frac{3-\sqrt{5}}{2},\frac{\sqrt{5}-1}{2}\}$, then there exists a $\mathcal{E}_{1}^{2}$-sentence $\varphi$ such that ${\sf P}(G(n,p)\models\varphi)$ does not converge. This finishes the proof of Theorem~\ref{thm_1_2}.\\

{\it Remark.} For every $c\in(0,1)$, there exists $p=\frac{3-\sqrt{5}}{2}+o(1)$ and a $\mathcal{E}_{1}^{2}$-sentence $\varphi$ such that the limit $\lim_{n\to\infty}{\sf P}(G(n,p)\models\varphi)$ equals $c$. Indeed, let $\varphi=\varphi_C^3$. From Lemma~\ref{main_lem} and the fact that $\varphi$ expresses the decreasing property, it follows that for every $\varepsilon$ and $n$ large enough there exists $p_{\varepsilon}(n)\in[\frac{3-\sqrt{5}}{2}-\varepsilon,\frac{3-\sqrt{5}}{2}+\varepsilon]$ such that ${\sf P}(G(n,p_{\varepsilon})\models\varphi)=c$. Therefore, there exists $p=\frac{3-\sqrt{5}}{2}+o(1)$ such that $\lim_{n\to\infty}{\sf P}(G(n,p)\models\varphi)=c$. The same results hold for some $p=\frac{\sqrt{5}-1}{2}+o(1)$.

\section{A sentence with 1 monadic variable and a sentence with 3 first order variables}
\label{examples}

Consider two rooted trees $F_1$ and $F_2$ with roots $R_1$ and $R_2$ respectively. Let us define {\it the product} of the rooted trees $F_1\cdot F_2$ in the following way. Let $E$ be the set of all possible edges $\{u,v\}$ where $u\in V(F_1)$, $v\in V(F_2)$ and $u$, $v$ are at the same distance from $R_1$, $R_2$ in $F_1$, $F_2$ respectively. Then $F_1\cdot F_2$ is the graph with the set of vertices $V(F_1)\sqcup V(F_2)$ and the set of edges $E(F_1)\sqcup E(F_2)\sqcup E$. Fix an arbitrary positive integer $a$ and consider two trees $F_1,F_2$ with $a$ and $\frac{4^a-1}{3}$ vertices respectively: $F_1$ is a simple path rooted at one of its end-points, and $F_2$ is a perfect $4$-ary tree (every non-leaf vertex of $F_2$ has 4 children and, for every $i\in\{1,\ldots,a-1\}$, the number of vertices at the distance $i$ from $R$ equals $4^i$). Denote $W_a=F_1\cdot F_2$.

Obviously, $v(W_a)=a+\frac{4^a-1}{3}$, and $e(W_a)=a+2\frac{4^a-1}{3}-2$.

\begin{lemma}
Consider an increasing sequence of positive integers $n_i$. Denote
$$
 k_i=\frac{2}{\ln[1/(1-p)]}\ln n_i-\frac{2}{\ln[1/(1-p)]}\ln\left[\frac{\sqrt[4]{24}\sqrt{(1-p)^5}}{p^2}\right].
$$
\begin{enumerate}
\item Let $0<c<1/2$, $\varepsilon>0$. If, for $i$ large enough, there is no integer $a_i$ such that $a_i+\frac{4^{a_i}-1}{3}\in(c k_i,k_i+\varepsilon)$, then a.a.s., for every $a$, there is no induced copy $F$ of $W_a$ in $G(n_i,p)$ such that every vertex outside $F$ has a neighbor inside $F$.

\item Let $1/2<C_1<C_2<1$. If, for $i$ large enough, there exists an integer $a_i$ such that $C_1 k_i\leq a_i+\frac{4^{a_i}-1}{3}\leq C_2 k_i$, then a.a.s. there is an induced copy $F$ of $W_{a_i}$ in $G(n_i,p)$ such that every vertex outside $F$ has a neighbor inside $F$. 
\end{enumerate}
\label{maximal_paths_mod}
\end{lemma}

{\it Proof of Lemma~\ref{maximal_paths_mod}.} 1. Let $W(a)$ be the number of induced copies of $W_a$ in $G(n,p)$. Let $s=a+\frac{4^a-1}{3}$. For $s=O(\ln n)$,
$$
 {\sf E}W(a)={n\choose s}\frac{s!}{24^{(s-a-1)/4}}p^{2s-a-2}(1-p)^{{s\choose 2}-2s+a+2}=
$$
$$
 e^{s\ln n-s\ln\left[\frac{\sqrt[4]{24}\sqrt{(1-p)^5}}{p^2}\right]-\frac{s^2}{2}\ln(1/(1-p))+O(\ln\ln n)}.
$$
Fix $\varepsilon>0$. For every $i$, let $a_i$ be the minimum number such that $s_i=a_i+\frac{4^{a_i}-1}{3}\geq\lceil k_i+\varepsilon\rceil$. Then ${\sf E}W(a_i)\leq e^{-\varepsilon\ln n_i+O(\ln\ln n_i)}$. Therefore, ${\sf P}(W(a_i)>0)\leq{\sf E}W(a_i)\to 0$ as $i\to\infty$.

So, it is enough to prove that a.a.s., for every set $X$ on at most $ck_i$ vertices, there is a vertex outside $X$ which has no neighbors inside $X$. The probability of this event is at least
$$
 1-{n_i\choose \lfloor ck_i\rfloor}(1-(1-p)^{\lfloor ck_i\rfloor})^{n_i-\lfloor ck_i\rfloor}\geq 1-e^{-A n_i^{1-2c}
 }\to 1\text{ as }n\to\infty
$$
for some positive constant $A$.\\

2. Now, let $C_1 k_i\leq s_i=a_i+\frac{4^{a_i}-1}{3}\leq C_2 k_i$. In what follows, we write $s,a,n$ instead of $s_i,a_i,n_i$ respectively. Let $\tilde W(a)$ be the number of induced copies of $W_a$ in $G(n,p)$ such that every vertex outside a copy has a neighbor inside. Then
$$
 {\sf E}\tilde W(a)={n\choose s}\frac{s!}{24^{(s-a-1)/4}}p^{2s-a-2}(1-p)^{{s\choose 2}-2s+a+2}(1-(1-p)^s)^{n-s}\geq
$$
$$
 e^{2(C_1-C_1^2)\log_{1/(1-p)}n\ln n+O(\ln n)}\to\infty.
$$

It remains to prove that $\frac{{\sf D}\tilde W(a)}{({\sf E}\tilde W(a))^2}\to 0$. We get
$$
 {\sf D}\tilde W(a)\leq{\sf E}\tilde W(a)+
 \sum_{\ell=1}^{s-1}F_{\ell}+o([{\sf E}\tilde W(a)]^2),
$$
where
$$
 F_{\ell}={n\choose s}{s\choose\ell}{n-s\choose s-\ell}\left[\frac{s!}{24^{(s-a-1)/4}}\right]^2\left(p^{2s-a-2}(1-p)^{{s\choose 2}-2s+a+2}\right)^2(1-p)^{-{\ell\choose 2}}G_{\ell},
$$
$$
G_{\ell}=\max\left\{1,\left[\frac{1-p}{p}\right]^{2\ell}\right\}.
$$
For every $\ell\in\{1,\ldots,s-1\}$,
$$
 \frac{F_{\ell}}{({\sf E}\tilde W(a))^2}=\frac{{s\choose \ell}{n-s\choose s-\ell}(1-p)^{-{\ell\choose 2}}G_{\ell}}{{n\choose s}}\leq
 \frac{(s/\ell)^{\ell}e^{\ell}{n\choose s-\ell}(1-p)^{-{\ell\choose 2}}G_{\ell}}{{n\choose s}}\leq
$$
$$
 \left(\frac{s^2\ell}{n}(1-p)^{-(\ell-1)/2}G_{\ell}\right)^{\ell(1+o(1))}\leq
 n^{(C_2-1)\ell(1+o(1))}.
$$
Therefore,
$$
\frac{{\sf D}\tilde W(a)}{({\sf E}\tilde W(a))^2}\leq\frac{1}{{\sf E}\tilde W(a)}+
sn^{(C_2-1)\ell(1+o(1))}+o(1)\to 0.\quad\Box
$$

Let us finish the proof of the first part of Theorem~\ref{thm_1_and_3}. We want to construct a sentence $\varphi_1$ which a.a.s. says that there is an induced copy of $W_a$ such that every vertex outside this copy has a neighbor in it. Note that such a copy has $a+\frac{4^a-1}{3}$ vertices. This is why we {\it do not need infinite number of disjunctions} for doing that! 

Let $\varphi_1=\exists X\,\phi(X)$, where $\phi=\phi(X)$ is a first order sentence with two binary predicates $\sim,=$ and one unary predicate $X$ saying that, for some $a$, the induced subgraph on $[X]:=\{v:\,X(v)\}$ is isomorphic to $W_a$ and every vertex outside $[X]$ has a neighbor inside $[X]$. It can be written, for example, in the following way:
$$
 \phi=\exists x\exists y^1\exists y_1^2\ldots\exists y_4^2\exists z_1\exists z_2\exists w\exists h\quad
 \mathrm{DEG}(x,\ldots,h)\wedge\mathrm{PATH}(z_2,w)\wedge\mathrm{PROD}(y_1,y^2_1,\ldots,y^2_4)\wedge
$$
$$
 \mathrm{PERFECT}(z_1,z_2,y_1^2,\ldots,y_4^2)\wedge\mathrm{TREE}(x,z_1)\wedge\mathrm{START}(z_1,z_2)
 \wedge\mathrm{MAX},
$$
where $x$ is the end-point (root) of the simple path $F_1$ which is adjacent to the root $y^1$ of $F_2$; $y_1^2,\ldots,y_4^2$ are children of $y^1$ in $F_2$; $z_1$ is the child of $x$ in $F_1$; $z_2$ is the child of $z_1$ in $F_1$; $w$ is the second end-point of $F_1$, and $h$ is its parent. Everywhere below, we skip predicates $X(\ell)$ for all first order variables $\ell$.

$\mathrm{DEG}(x,\ldots,h)$ defines degrees of all the vertices in the induced subgraph on the set $[X]$ and the relations between the distinguished vertices:
$$
 \mathrm{DEG}(x,\ldots,h)=(d(x)=2)\wedge(d(y^1)=5)\wedge\left[\bigwedge_{\ell\in\{z_2,w,h\}}(d(\ell)\geq 7)\right]\wedge\left[\bigwedge_{\ell\in\{z_1,y^2_1,\ldots,y^2_4\}}(d(\ell)=6)\right]\wedge
$$
$$
 (\forall u\,\,[u\sim w\wedge u\neq h]\rightarrow[d(u)=2])\wedge
 (\forall u\,\,[u\nsim w\wedge u\neq x\wedge u\neq y^1]\rightarrow[d(u)\geq 6])\wedge
$$
$$
 (x\sim y^1)\wedge(x\sim z_1)\wedge (z_1\sim z_2)\wedge(h\sim w)\wedge\left[\bigwedge_{i=1}^4(z_1\sim y_i^2)\wedge(y^1\sim y_i^2)\right]\wedge\left[\bigwedge_{1\leq i<j\leq 4}(y_i^2\neq y_j^2)\right].
$$
Here, the first order sentence $d(\ell)=m$ ($d(\ell)\geq m$) says that the degree of $\ell$ in the induced subgraph on the set $[X]$ equals $m$ (at least $m$).

$\mathrm{START}(z_1,z_2)$ says that each of the four children of the root in $F_2$ has four children in $F_2$, and each of them is adjacent to $z_2$:
$$
 \mathrm{START}(z_1,z_2)=\forall u\,\,\biggl[u\sim z_1\wedge d(u)=6\biggr]\rightarrow
$$
$$ 
 \left[\exists v_1\ldots\exists v_4\,\,\left(\bigwedge_{i\neq j\in\{1,\ldots,4\}}[v_i\sim z_2\wedge v_i\neq z_1\wedge v_i\neq v_j]\right)\right].
$$

$\mathrm{PATH}(z_2,w)$ says that $F_1$ is either a simple path or a union of a simple path and simple cycles:
$$
 \mathrm{PATH}(z_2,w)=
 \biggl(\bigwedge_{\ell\in\{z_2,w\}}[\exists u\,\,(u\sim\ell)\wedge(d(u)\geq 7)\wedge(\forall v\,\,[d(v)\geq 7\wedge v\sim\ell]\rightarrow[v=u])]\biggr)\wedge
$$
$$
 \biggl(\forall u\quad[(d(u)\geq 7)\wedge(u\neq w)\wedge(u\neq z_2)]\rightarrow
$$
$$
 \biggl[\exists u_1\exists u_2\,\,(u_1\neq u_2)\wedge\left(\bigwedge_{i=1,2} [u_i\sim u\wedge d(u_i)\geq 7]\right)\wedge
$$
$$ 
 \left(\forall v\,\,[d(v)\geq 7\wedge v\sim u]\rightarrow\left[\bigvee_{i=1,2}v=u_i\right]\right)\biggr]\biggr).
$$
$\mathrm{PROD}(y_1,y^2_1,\ldots,y^2_4)$ says that every vertex of $F_2$ has an only neighbor in $V(F_1)$:
$$
 \mathrm{PROD}(y_1,y^2_1,\ldots,y^2_4)=\left(\bigwedge_{i=1}^4[\forall u\quad (u\sim y^2_i\wedge u\neq y^1)\rightarrow(d(u)=6)]\right)\wedge
$$
$$
 \biggl(\forall u\,\,\left[(d(u)=6)\wedge\left(\bigwedge_{i=1}^4
 u\neq y^2_i\right)\right]\rightarrow 
$$
$$ 
 [\exists v\,\,(v\sim u)\wedge (d(v)\geq 7)\wedge(\forall \tilde v\,\,[d(\tilde v)\geq 7\wedge v\sim u]\rightarrow[v=\tilde v])]\biggr).
$$
$\mathrm{PERFECT}(z_1,z_2,y_1^2,\ldots,y_4^2)$ is the main formula that defines the structure of the graph $F_1\cdot F_2$. It says that every vertex $u$ of $F_2$ (except for the root and the leaves) has a neighbor in $F_1$ that have two $v_1$ and $v_2$ neighbors in $F_1$ such that $v_1$ has an only common neighbor with $u$ in $F_2$, and $v_2$ has four common neighbors with $u$ in $F_2$:
$$
 \mathrm{PERFECT}(z_1,z_2,y_1^2,\ldots,y_4^2)=\forall u\quad \left[(d(u)=6)\wedge(u\neq z_1)\wedge\left(\bigwedge_{i=1}^4 u\neq y^2_i\right)\right]\rightarrow
$$
$$
 \biggl[\exists v\,\,(d(v)\geq 7)\wedge(v\sim u)\wedge
$$
$$
 \biggl(\exists v_1\,\,[(v_1=z_1 \wedge v=z_2)\vee (d(v_1)\geq 7\wedge v\neq z_2)]\wedge[N_6(v_1,u)=1]\wedge[N(v_1,u)=2]\biggr)\wedge
$$
$$
 \biggl(\exists v_2\,\,[d(v_2)\geq 7]\wedge[N(v_2,u)=5]\wedge[(v_2=w\wedge N_2(v_2,u)=4)\vee(v_2\neq w\wedge N_6(v_2,u)=4)]\biggr)
 \biggr].
$$
Here, the first order sentences $N(\ell_1,\ell_2)=n$ and $N_m(\ell_1,\ell_2)=n$ say that $\ell_1,\ell_2$ have exactly $n$ common neighbors and exactly $n$ common neighbors with the degree $m$ respectively.

$\mathrm{TREE}(x,z_1)$ says that $F_2$ is a tree (two vertices at the same distance from the root are not adjacent, and have only one common neighbor in $F_2$) whenever $\mathrm{DEG}(x,\ldots,h)$ and $\mathrm{START}(z_1,z_2)$ are true:
$$
 \mathrm{TREE}(x,z_1)=
 \forall u_1\forall u_2\quad\left[\left(\bigwedge_{i=1,2}([d(u_i)=6 \wedge u_i\neq z_1]\vee[d(u_i)=2\wedge u_i\neq x])\right)\wedge\right.
$$
$$
 \left.\left(\exists v\,\, [v=z_1\vee d(v)\geq 7]\wedge\left[\bigwedge_{i=1,2}u_i\sim v\right]\right)\right]\rightarrow\biggl[(u_1\nsim u_2)\wedge(N(u_1,u_2)=2)\biggr].
$$
Note that $\mathrm{TREE}(x,z_1)$, $\mathrm{PERFECT}(z_1,z_2,y_1^2,\ldots,y_4^2)$ and $\mathrm{PROD}(y_1,y^2_1,\ldots,y^2_4)$ imply that $F_1$ is a simple path whenever $\mathrm{PATH}(z_2,w)$ is true.

$\mathrm{MAX}$ says that every vertex outside $[X]$ has a neighbor inside $[X]$.

It remains to consider two sequences $n_i$, $m_i$ such that $\lim_{i\to\infty}{\sf P}(G(n_i,p)\models\varphi_1)=1$, and $\lim_{i\to\infty}{\sf P}(G(m_i,p)\models\varphi_1)=0$. By Lemma~\ref{maximal_paths_mod},
$$
n_i=\left\lceil \frac{\sqrt[4]{24}\sqrt{(1-p)^5}}{p^2} \left(\frac{1}{1-p}\right)^{\frac{3i+4^i-1}{4}}\right\rceil,\quad
m_i=\left\lceil \frac{\sqrt[4]{24}\sqrt{(1-p)^5}}{p^2} \left(\frac{1}{1-p}\right)^{\frac{3i+4^i-1}{2}}\right\rceil
$$
are the required sequences.\\

It remains to prove that the same can be done using an EMSO sentence $\varphi_2$ with arbitrary number of existential monadic variables and only three first order variables. But this is much easier. Let $\varphi_2=\exists X\exists A\exists B\exists W\left[\bigwedge_{i=1}^3\exists P_i\left(\bigwedge_{j=1}^4\exists C_i^j\quad\phi\right)\right]$, where $\phi$ is a first order sentence with two binary predicates $\sim,=$ and unary predicates $X,A,B,W,P_i,C_i^j$ saying that (everywhere below, for $i\in\{1,2,3\}$, the values $i-1$ and $i+1$ are still in $\mathbb{Z}_3+1$: $1-1=3$, $3+1=1$)

\begin{enumerate}

\item every vertex outside $[X]$ has a neighbor inside $[X]$;

\item the sets $[A], [B], [W], [P_i], [C_i^j]$ form a partition of $[X]$;

\item each of the sets $[A]$, $[B]$, $[W]$ contains only one vertex;

\item the only vertex of $[A]$ is adjacent to the only vertex of $[B]$;

\item the only vertex of $[A]$ is not adjacent to any vertex of $[W],[P_2],[P_3],[C_i^j]$;

\item the only vertex of $[B]$ is not adjacent to any vertex of $[W],[P_i],[C_2^j],[C_3^j]$;

\item every vertex in $[C_i^j]$ has exactly one neighbor in $[W]\sqcup[P_1]\sqcup[P_2]\sqcup[P_3]$, and this neighbor belongs either to $[W]$ or to $[P_i]$;

\item for every $j\in\{1,2,3,4\}$, if a vertex $c$ from $[C_1^j]$ is adjacent to the only vertex in $[B]$, then 

\begin{itemize}
\item $c$ has a neighbor in $[P_1]$ which is adjacent to the only vertex from $[A]$, and, 

\item for every $\tilde j\in\{1,2,3,4\}$, $c$ has a neighbor in $C_2^{\tilde j}$ which is adjacent to a vertex from $[P_2]$, which, in turn, is adjacent to a neighbor of $c$ in $[P_1]$;
\end{itemize}

\item for every $j\in\{1,2,3,4\}$ and every $i\in\{1,2,3\}$, if a vertex $c$ from $[C_i^j]$ is adjacent to the only vertex in $[W]$, then, 

\begin{itemize}
\item for some $\tilde j\in\{1,2,3,4\}$, $c$ has a neighbor in $[C_{i-1}^{\tilde j}]$ which, in turn, has a neighbor in $[P_{i-1}]$, which is adjacent to the only vertex in $[W]$;
\end{itemize}

\item for every $j\in\{1,2,3,4\}$ and every $i\in\{1,2,3\}$, if a vertex $c$ from $[C_i^j]$ is adjacent to a vertex $p$ in $[P_i]$, $c$ is not adjacent to the only vertex in $[B]$, and $c$ is not adjacent to a neighbor of the only vertex in $[W]$, then

\begin{itemize}
\item $p$ is adjacent to a vertex from $[P_{i-1}]$ which, in turn, for some $\tilde j\in\{1,2,3,4\}$, is adjacent to a vertex from $[C_{i-1}^{\tilde j}]$, which is a neighbor of $c$, and, 

\item for every $\tilde j\in\{1,2,3,4\}$, $c$ has a neighbor in $C_{i+1}^{\tilde j}$ which is adjacent to a vertex from $[P_{i+1}]$, which, in turn, is adjacent to $p$;
\end{itemize}

\item the graph induced on $[A]\sqcup[W]\sqcup[P_1]\sqcup[P_2]\sqcup[P_3]$ is a simple path with the end-points in $[A]$ and $[W]$:

\begin{enumerate}

\item the only vertex of $[A]$ is adjacent to a vertex from $[P_1]$, and there is only one such vertex in $[P_1]$,


\item there exists $i\in\{1,2,3\}$ such that the only vertex $w$ of $[W]$ is adjacent to a vertex from $[P_i]$, and there is no other neighbors of $w$ in $[P_1]\sqcup[P_2]\sqcup[P_3]$,

\item if a vertex $p$ from $[P_1]$ is adjacent to the only vertex in $[A]$, then it is also adjacent to a vertex from $[P_2]$, there is only one such vertex in $[P_2]$, and every vertex from $[P_3]$ is not adjacent to $p$,

\item for each $i\in\{1,2,3\}$, if a vertex $p$ from $[P_i]$ is adjacent to the only vertex in $[W]$, then it is also adjacent to a vertex from $[P_{i-1}]$, there is only one such vertex in $[P_{i-1}]$, and every vertex from $[P_{i+1}]$ is not adjacent to $p$,

\item for each $i\in\{1,2,3\}$, if a vertex $p$ from $[P_i]$ is not adjacent to the only vertex in $[A]$ and is not adjacent to the only vertex in $[W]$, then, for each $\tilde i\in\{1,2,3\}\setminus\{i\}$, $p$ is adjacent to a vertex from $[P_{\tilde i}]$, there is only one such vertex in $[P_{\tilde i}]$,

\item for each $i\in\{1,2,3\}$, $[P_i]$ is an independent set;

\end{enumerate}

\item the graph induced on $[B]\sqcup\bigcup_{i,j}[C_i^j]$ is a perfect $4$-ary tree with the root in $[B]$:

\begin{enumerate}

\item for every $j\in\{1,2,3,4\}$, the only vertex of $[B]$ is adjacent to a vertex from $[C_1^j]$, and there is only one such vertex in $[C_1^j]$,

\item for every $j\in\{1,2,3,4\}$, if a vertex $c$ from $[C_1^j]$ is adjacent to the only vertex in $[B]$, then, for every $\tilde j\in\{1,2,3,4\}$, it is also adjacent to a vertex from $[C_2^{\tilde j}]$, there is only one such vertex in $[C_2^{\tilde j}]$, and every vertex from $[C_3^{\tilde j}]$ is not adjacent to $c$,

\item for every $j\in\{1,2,3,4\}$ and every $i\in\{1,2,3\}$, if a vertex $c$ from $[C_i^j]$ is adjacent to the only vertex in $[W]$, then, for some $\tilde j\in\{1,2,3,4\}$, $c$ has a neighbor in $[C_{i-1}^{\tilde j}]$, and $c$ has no more neighbors in all the $C$-sets,

\item for every $j\in\{1,2,3,4\}$ and every $i\in\{1,2,3\}$, if a vertex $c$ from $[C_i^j]$ is not adjacent to the only vertex in $[B]$ and is not adjacent to the only vertex in $[W]$, then 

\begin{itemize}
\item there exist $\tilde j\in\{1,2,3,4\}$ such that $c$ has a neighbor in $[C_{i-1}^{\tilde j}]$, and $c$ has no more neighbors in all the $C_{i-1}$-sets,

\item for every $\tilde j\in\{1,2,3,4\}$, $c$ has exactly one neighbor in $[C_{i+1}^{\tilde j}]$,
\end{itemize}

\item for every $j\in\{1,2,3,4\}$ and every $i\in\{1,2,3\}$, $[C_i^j]$ is an independent set.

\end{enumerate}

\end{enumerate}

All the properties above can be easily expressed using 3 variables, and $\phi$ is the conjunction of them.

\section{Acknowledgements} This paper is supported by the grant of the Russian Science Foundation No. 16-11-10014. 

\end{document}